\documentclass[12pt]{article}
\usepackage{ams}
\title{METHOD LAC}
\author{Joachim Nzotungicimpaye\\Kigali Institute of Education, Department of Mathematics\\ P.O.Box 5039,Kigali-Rwanda\\e-mail kimpaye @kie.ac.rw}

\begin{document}
\maketitle
\begin{abstract}
This is a pedagogical paper which present a mnemotechnical method
 that we call LAC for Lists , Arrangements and Combinations. It can help students or any one to recollect formulae from combinatorial theory
(\cite{schaum}),\cite{boas}, \cite{barnett}, \cite{any})without an
a priori memorization of them.
\end{abstract}
\section{Introduction}
Let us suppose that we have a set $X$ of $n$ objects. The
following questions are generally addressed. In how many ways can
we organize them $p$ by $p$ if repetitions and order are accepted,
if repetitions are not accepted but order is still accepted and
finally if repetitions nor order are accepted? The possibilities
are respectively called \textbf{L}ists, \textbf{A}rrangements and
\textbf{C}ombinations. At secondary school students learn that
there are $n^p$ lists, $\frac{n!}{(n-p)!}$ arrangements and
$\frac{n!}{(n-p)!p!}$ combinations (\cite{schaum}),\cite{boas},
\cite{barnett}, \cite{any}). We give in the following section by
use of the theory of matrices a mnemo-technical method which
permit anyone to found out these numbers rapidly without having
memorized them. Using the three initials from Lists, Arrangements
and Combinations ,we propose that this method be called method
\textbf{LAC}
.\\
\\

\section{Method LAC}
We suppose that we have at our disposition $n$ objects to organize
two by two and find out the corresponding formulae by use of
matrices. We then generalize the results to the organizations $p$
by $p$.
\subsection{\textbf{L}ists}
When lists are admitted, it means that we accept repetition and
order, this situation is represented by a general $n$ by $n$
matrix . If we limit ourself to a set $X=\{a,b,c,d,e\} $ of five
objects ,we will have a matrix of the form
\begin{eqnarray}
g=\left (
\begin{array}{ccccc}
aa&ab&ac&ad&ae\\ba&bb&bc&bd&bd\\ca&cb&cc&cd&ce\\da&db&dc&dd&de\\ea&eb&ec&ed&ee
\end{array}
\right )
\end{eqnarray}
We then have $5^2$ possibilities and $n^2$ possibilities if we $n$
objects are considered. If we generalize this to the organization
$p$ by $p$ , we will then have $n^p$ lists.
\subsection{\textbf{A}rrangements and Permutations}
When arrangements are admitted, it means that we do not accept
repetition but we still accept order, this situation is
represented by a general $n$ by $n$ matrix where we have filled
out diagonal entries which represents $n$ non possibilities
because there are no repetitions. Limiting to five objects give
rise to a matrix like
\begin{eqnarray}
g=\left (
\begin{array}{ccccc}
~&ab&ac&ad&ae\\ba&~&bc&bd&bd\\ca&cb&~&cd&ce\\da&db&dc&~&de\\ea&eb&ec&ed&~
\end{array}
\right )
\end{eqnarray}
We then have $5^2-5=5(5-1)$ possibilities and $n(n-1)$
possibilities for $n$ objects in consideration. Let us recall that
\begin{eqnarray}
n(n-1)=\frac{n!}{(n-2)!}
\end{eqnarray}
where
\begin{eqnarray}
n!=n.(n-1).(n-2)......4.3.2.1
\end{eqnarray}
 By generalizing to the organization $p$ by $p$ , we obtain
$\frac{n!}{(n-p)!}$ arrangements. Note that a permutation of $n$
objects is an arrangement of $n$ objects $n$ by $n$. Using the
convention telling us that $0!=1$ , we see that there are $n!$
permutations of $n$ objects.

\subsection{\textbf{C}ombinations}
When combinations are admitted, it means that we do not accept
repetition nor order.This situation is represented by a lower or
upper $n$ by $n$ matrix where we have filled out the diagonal
entries and the upper entries or the lower entries because there
are no repetitions nor order. Limiting to five objects give rise
to a matrix like
\begin{eqnarray}
g=\left (
\begin{array}{ccccc}
~&ab&ac&ad&ae\\~&~&bc&bd&bd\\~&~&~&cd&ce\\~&~&~&~&de\\~&~&~&~&~
\end{array}
\right )
\end{eqnarray}
We then have half of the previous possibilities
,$\frac{5(5-1)}{2}$ possibilities, it means $\frac{n!}{(n-2)!2!}$
when $n$ objects are considered. The generalization to the
organization $p$ by $p$ give rise to $\frac{n!}{(n-p)!p!}$
combinations.


\begin{thebibliography}{99}
\bibitem{schaum}Seymour Lipschutz, {\it Theory and Problems of Finite Mathematics } , Schaum's Outline Series , McGraw-Hill Book Company,NY $1966$
\bibitem{boas} M.L.Boas, {\it Mathematical Methods in the Physical Sciences}, $2^{nd}$
Edition, John Wiley and Sons, $1983$.
\bibitem{barnett}R.A.Barnett and al., {\it College algebra, a
graphing approach}, McGraw-Hill, $2000$
\bibitem{any} Any other book on Combinatorial algebra.
\end{thebibliography}
\end{document}